\documentclass[journal,article,submit,moreauthors]{mdpi} 
\UseRawInputEncoding
\firstpage{1} 
\hreflink{https://doi.org/} 

\usepackage{wrapfig}
\usepackage[pscoord]{eso-pic}
\usepackage[fulladjust]{marginnote}

\reversemarginpar

\definecolor{Gray}{gray}{.25}

\Title{Eliminating oscillation in partial sum approximation of periodic function}

\Author{Shi-Lin Li $^{1}$, Yuan-Yuan Liu $^{1,3}$ , Wen-Du Li $^{2,3,\ddagger}$and Wu-Sheng Dai $^{1,\ddagger}$}

\address{%
$^{1}$ \quad Department of Physics, Tianjin University, Tianjin 300350, P.R. China\\
$^{2}$ \quad College of Physics and Materials Science, Tianjin Normal University, Tianjin 300387, PR China\\
$^{3}$ \quad Theoretical Physics Division, Chern Institute of Mathematics, Nankai University, Tianjin, 300071, P. R. China
}

\firstnote{liwendu@tjnu.edu.cn.} 
\secondnote{daiwusheng@tju.edu.cn.}

\abstract{If we cannot obtain all terms of a series, or if we cannot sum up a series, we
have to turn to the partial sum approximation which approximate a function by
the first several terms of the series. However, the partial sum approximation
often does not work well for periodic functions. In the partial sum
approximation of a periodic function, there exists an incorrect oscillation
which cannot be eliminated by keeping more terms, especially at the domain
endpoints. A famous example is the Gibbs phenomenon in the Fourier expansion.
In the paper, we suggest an approach for eliminating such oscillations in the
partial sum approximation of periodic functions.}


\begin{document}
\nolinenumbers



\section{Introduction}

If a series cannot be exactly summed up, one turns to approximate the series
by a partial sum, i.e., approximate the series by the sum of its first several
terms. In many cases, the accuracy of the partial sum increases as the number
of terms increases. But in the expansion of periodic functions, one encounters
such a situation: the accuracy cannot be improved by increasing the number of
terms in the partial sum. In other words, even if the number of the term of
the partial sum is increased, the accuracy will not be improved. In the
partial sum approximation of a periodic function, there exists an incorrect
oscillation, especially at the domain endpoints. An important example of such
kinds of problem is the Gibbs phenomenon of the Fourier expansion. In this
appendix, we suggest an approach to solve the problem. We will take the
problem encountered in the calculation of scattering cross section in the main
text as an example to illustrate this approach.

In section \ref{ex}, we illustrate the problem by an example. In section
\ref{Scheme}, we suggest a scheme for eliminating the oscillation. In sections
\ref{QM} and \ref{RN}, we consider examples in scattering. The conclusions are
summarized in section \ref{Conclusion}.

\section{Oscillation in partial sum approximation of periodic function
\label{ex}}

In order to illustrate the problem encountered in the expansion of periodic
functions, we first take a look at an example in which the sum can be exactly performed.

The periodic function
\begin{equation}
f\left(  \theta\right)  =\frac{1}{2\sin\frac{\theta}{2}}%
\end{equation}
can be expanded as%
\begin{equation}
f\left(  \theta\right)  =\sum_{l=0}^{\infty}P_{l}\left(  \cos\theta\right)
,\label{sigamp}%
\end{equation}
where $P_{l}\left(  \cos\theta\right)  $ is the Legendre polynomial. If we use
an $N$-term partial sum to approximate the exact result,%
\begin{equation}
f_{N}\left(  \theta\right)  =\sum_{l=0}^{N}P_{l}\left(  \cos\theta\right)
,\label{fNsin}%
\end{equation}

\vspace*{-0.3cm}
\begin{wrapfigure}[18]{l}[0.5cm]{0.8\textwidth}
\includegraphics[width=0.8\textwidth]{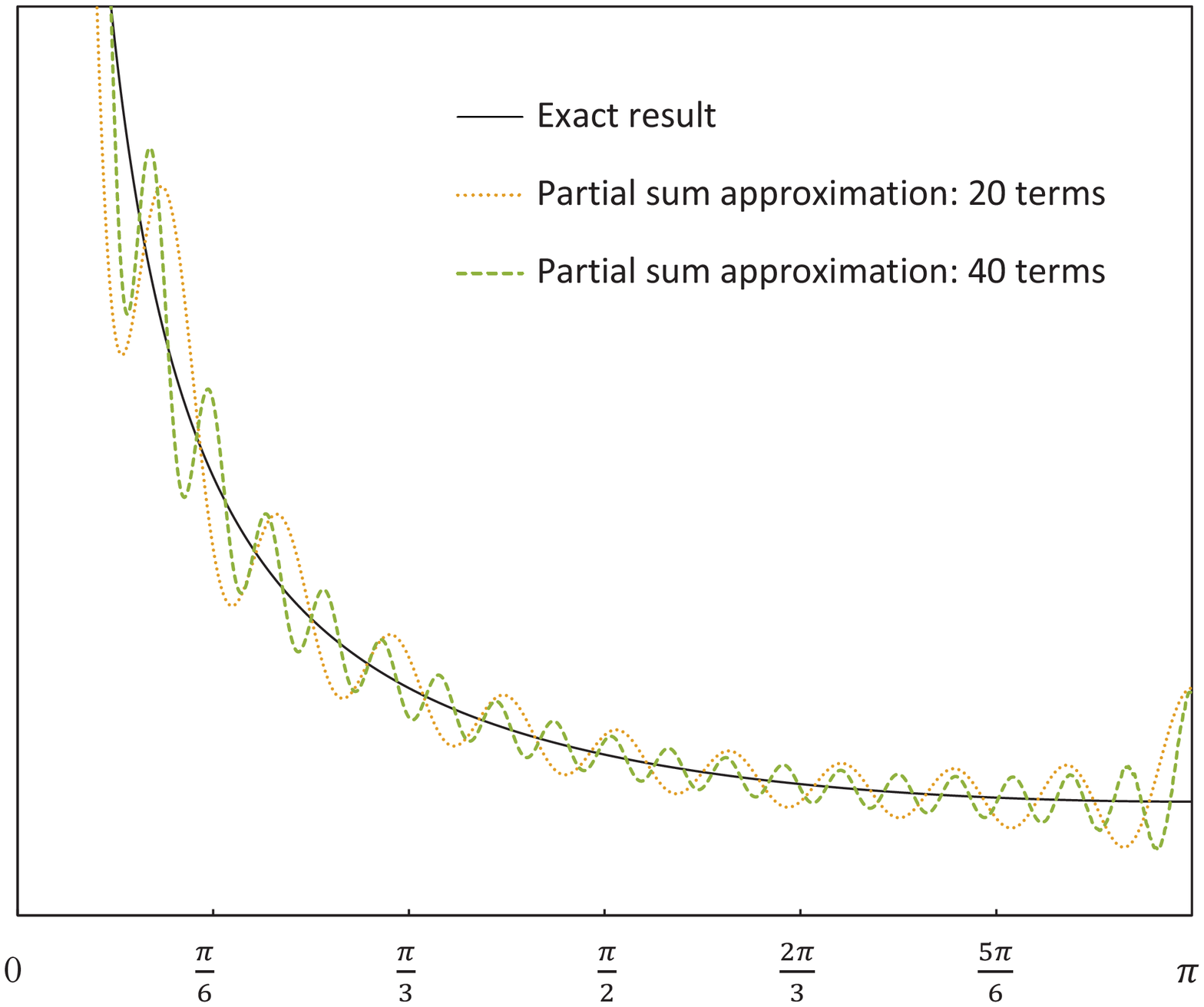}
\captionsetup{labelformat=empty} 
\caption{} 
\label{sum.eps} 
\end{wrapfigure}
\
\qquad\qquad\\
\qquad\qquad\\
\qquad\qquad\\
\qquad\qquad\\
\qquad\qquad\\
\qquad\qquad\\
\qquad\qquad\\
\qquad\qquad\\
\qquad\qquad\\
\vspace*{4.5cm}
\qquad\qquad\\
\qquad\qquad\\
\qquad\qquad\\
\qquad\qquad\\
\qquad\qquad\\
\qquad\qquad\\
\qquad\qquad\\

\bigskip
\noindent an incorrect oscillation appears Fig. \ref{sum.eps}. This oscillation does not
exist in the exact result, and can not be eliminated by increasing the number
of partial sum terms, as shown in Fig. \ref{sum.eps}. 
Especially at the edges, i.e., as $\theta$ approaches $\pi$, the
error does not decrease as the total number of terms of the partial sum
increases. This oscillation is essentially the Gibbs phenomenon in the Fourier expansion.
\switchcolumn
\newpage
\vspace*{4cm}
{\color{Gray} 
\textbf{Figure \ref{sum.eps}.} The incorrect
oscillation in partial sum approximation of the periodic function $\displaystyle f\left(
\theta\right)  =\frac{1}{2\sin\frac{\theta}{2}}$: The exact result does not
oscillate, but the partial sum approximation does. 
}
\switchcolumn

\section{Generalized Pad{\'{e} }approximant \label{Scheme}}

The example in section \ref{ex} is a special case of the problem encountered
in scattering theory. Generally, in scattering theory, when using the partial
wave method to calculate the scattering cross section, we encounter the
following sum:%
\begin{equation}
f\left(  \theta\right)  =\sum_{l=0}^{\infty}c_{l}P_{l}\left(  \cos
\theta\right)  , \label{fP}%
\end{equation}
where $c_{l}$ is the expansion coefficient and $P_{l}\left(  z\right)  $ is
the Legendre polynomial.

If this sum cannot be performed exactly, we have to truncate the series and
approximate the function $f\left(  \theta\right)  $ by a partial sum
consisting of the first $N$ terms of the series:%
\begin{equation}
f_{N}\left(  \theta\right)  =\sum_{l=0}^{N}c_{l}P_{l}\left(  \cos
\theta\right)  . \label{Nterms}%
\end{equation}

The basis of this expansion is the Legendre polynomial, while the basis of the
Fourier series is the sine/cosine function or the exponential function. The
Fourier series can be obtained from a power series $\left\{  x^{n}\right\}  $
by replacing the basis of the Fourier series $x$ with $e^{in\theta}$. The
Legendre polynomial is an orthogonalized power series, which is a rearranged
power series, so the series with $\left\{  P_{l}\left(  \cos\theta\right)
\right\}  $ as the basis is a rearranged Fourier series, which is essentially
still a Fourier expansion. Therefore, the problem encountered here is
essentially the Gibbs phenomenon in the Fourier expansion.

In this appendix, we will construct a modified Pad{\'{e}} approximant to solve
this problem.

The Pad{\'{e} }approximant is to use a rational function instead of the power
series to approximate a function \cite{tian2021pade}. Since the common
Pad{\'{e} }approximant is of low efficiency in this case, in this appendix, we
introduce a generalized Pad{\'{e} }approximant to approximate the series
(\ref{Nterms}).

The generalized Pad{\'{e} }approximant is constructed as
\begin{equation}
f_{\left[  L/M\right]  }\left(  \theta\right)  =\frac{\sum_{n=0}^{L}a_{n}%
P_{n}\left(  \cos\theta\right)  }{\sum_{m=0}^{M}b_{m}P_{m}\left(  \cos
\theta\right)  }. \label{PadeL}%
\end{equation}
The generalized Pad{\'{e} }approximant is to approximate the function
$f\left(  \theta\right)  $ with the rational form (\ref{PadeL}) instead of the
polynomial form (\ref{Nterms}). The numerator of the rational expression is a
polynomial of order $L$ and the denominator is a polynomial of order $M$. In
principle, $L$ and $M$ can be chosen arbitrarily as long as the condition
$L+M=N$ is satisfied. Unlike the power series expansion, there is no unified
method to determine the coefficients $a_{n}$ and $b_{m}$ in the rational
expression (\ref{PadeL}).

In the Pad{\'{e} }approximant, one uses the power series expansion to
determine the coefficients of the rational approximation. Similarly, we here
use the partial sum (\ref{Nterms}) to determine the coefficients of the
rational form approximation (\ref{PadeL}), i.e., to determine $a_{n}$ and
$b_{m}$ from $c_{l}$ by equaling the rational expression (\ref{PadeL}) and the
polynomial (\ref{Nterms}):%
\begin{equation}
\frac{\sum_{n=0}^{L}a_{n}P_{n}\left(  \cos\theta\right)  }{\sum_{m=0}^{M}%
b_{m}P_{m}\left(  \cos\theta\right)  }=\sum_{l=0}^{L+M}c_{l}P_{l}\left(
\cos\theta\right)  . \label{LPade}%
\end{equation}
These are equations determining the coefficients $a_{n}$ and $b_{m}$.

Next, we solve the coefficients $a_{n}$ and $b_{m}$ from Eq. (\ref{LPade}).

By Eq. (\ref{LPade}) we have%
\begin{align}
\sum_{n=0}^{L}a_{n}P_{n}\left(  \cos\theta\right)   &  =\left(  \sum
_{l=0}^{L+M}c_{l}P_{l}\left(  \cos\theta\right)  \right)  \left(  \sum
_{m=0}^{M}b_{m}P_{m}\left(  \cos\theta\right)  \right) \nonumber\\
&  =\sum_{l=0}^{L+M}\sum_{m=0}^{M}c_{l}b_{m}P_{l}\left(  \cos\theta\right)
P_{m}\left(  \cos\theta\right)  . \label{coeeq}%
\end{align}
The coefficients $a_{n}$ and $b_{m}$ are given by equaling the coefficients of
the Legendre polynomial of the same order on both sides.

The coefficient $a_{n}$ can be obtained by utilizing the orthogonality of the
Legendre polynomial by multiplying $P_{n}\left(  \cos\theta\right)  $ on both
sides of Eq. (\ref{coeeq}) and performing the integral $\int_{-1}^{1}%
d\cos\theta$:%
\begin{equation}
a_{n}=\sum_{l=0}^{L+M}\sum_{m=0}^{M}c_{l}b_{m}\int_{-1}^{1}d\cos\theta
P_{n}\left(  \cos\theta\right)  P_{l}\left(  \cos\theta\right)  P_{m}\left(
\cos\theta\right)  .
\end{equation}
Using%
\begin{equation}
\int_{-1}^{1}d\cos\theta P_{l}\left(  \cos\theta\right)  P_{m}\left(  \cos\theta\right)
P_{n}\left(  \cos\theta\right)  =2\left(
\begin{array}
[c]{ccc}%
l & m & n\\
0 & 0 & 0
\end{array}
\right)  ^{2}%
\end{equation}
with $\left(
\begin{array}
[c]{ccc}%
l & m & n\\
0 & 0 & 0
\end{array}
\right)  $ the $3j$ coefficient \cite{olver2010nist}, we have
\begin{equation}
a_{n}=2\sum_{l=0}^{L+M}\sum_{m=0}^{M}c_{l}b_{m}\left(
\begin{array}
[c]{ccc}%
l & m & n\\
0 & 0 & 0
\end{array}
\right)  ^{2}. \label{aneq}%
\end{equation}

From the right side of Eq. (\ref{coeeq}), it can be seen that the highest
order of the Legendre polynomial in the sum is $L$, which requires%
\[
a_{n}=0\text{ \ \ for }L+1\leq n\leq L+M.
\]
That is
\begin{equation}
2\sum_{l=0}^{L+M}\sum_{m=0}^{M}c_{l}b_{m}\left(
\begin{array}
[c]{ccc}%
l & m & n\\
0 & 0 & 0
\end{array}
\right)  ^{2}=0\text{ \ \ for }L+1\leq n\leq L+M. \label{an}%
\end{equation}
This is a system of linear equations that determines $b_{l}$:%
\end{paracol}
\begin{align}
b_{1}\sum_{m=0}^{L+M}c_{m}\left(
\begin{array}
[c]{ccc}%
1 & m & L+1\\
0 & 0 & 0
\end{array}
\right)  ^{2}+\cdots+b_{M}\sum_{m=0}^{L+M}c_{m}\left(
\begin{array}
[c]{ccc}%
M & m & L+1\\
0 & 0 & 0
\end{array}
\right)  ^{2}  &  =-b_{0}\sum_{m=0}^{L+M}c_{m}\left(
\begin{array}
[c]{ccc}%
0 & m & L+1\\
0 & 0 & 0
\end{array}
\right)  ^{2},\nonumber\\
b_{1}\sum_{m=0}^{L+M}c_{m}\left(
\begin{array}
[c]{ccc}%
1 & m & L+2\\
0 & 0 & 0
\end{array}
\right)  ^{2}+\cdots+b_{M}\sum_{m=0}^{L+M}c_{m}\left(
\begin{array}
[c]{ccc}%
M & m & L+2\\
0 & 0 & 0
\end{array}
\right)  ^{2}  &  =-b_{0}\sum_{m=0}^{L+M}c_{m}\left(
\begin{array}
[c]{ccc}%
0 & m & L+2\\
0 & 0 & 0
\end{array}
\right)  ^{2},\nonumber\\
&  \vdots\nonumber\\
b_{1}\sum_{m=0}^{L+M}c_{m}\left(
\begin{array}
[c]{ccc}%
1 & m & L+M\\
0 & 0 & 0
\end{array}
\right)  ^{2}+\cdots+b_{M}\sum_{m=0}^{L+M}c_{m}\left(
\begin{array}
[c]{ccc}%
M & m & L+M\\
0 & 0 & 0
\end{array}
\right)  ^{2}  &  =-b_{0}\sum_{m=0}^{L+M}c_{m}\left(
\begin{array}
[c]{ccc}%
0 & m & L+M\\
0 & 0 & 0
\end{array}
\right)  ^{2}. \label{beq}%
\end{align}
\begin{paracol}{2}
\switchcolumn
These $M$ equations solve $M$ coefficients $b_{l}$. For a rational expression,
the value of $b_{0}$ can be taken arbitrarily, so for convenience we take
$b_{0}=1$.

The solution of Eq. (\ref{beq}) is
\begin{equation}
b_{k}=\frac{\det B_{k}}{\det A}. \label{coefofb}%
\end{equation}
Here the matrix $A$ is the coefficient matrix of Eq. (\ref{beq}),%
\begin{equation}
A=\left(
\begin{array}
[c]{ccc}%
\sum\limits_{m=0}^{L+M}c_{m}\left(
\begin{array}
[c]{ccc}%
1 & m & L+1\\
0 & 0 & 0
\end{array}
\right)  ^{2} & \cdots & \sum\limits_{m=0}^{L+M}c_{m}\left(
\begin{array}
[c]{ccc}%
M & m & L+1\\
0 & 0 & 0
\end{array}
\right)  ^{2}\\
\sum\limits_{m=0}^{L+M}c_{m}\left(
\begin{array}
[c]{ccc}%
0 & m & L+2\\
0 & 0 & 0
\end{array}
\right)  ^{2} & \cdots & \sum\limits_{m=0}^{L+M}c_{m}\left(
\begin{array}
[c]{ccc}%
M & m & L+2\\
0 & 0 & 0
\end{array}
\right)  ^{2}\\
\vdots & \ddots & \vdots\\
\sum\limits_{m=0}^{L+M}c_{m}\left(
\begin{array}
[c]{ccc}%
1 & m & L+M\\
0 & 0 & 0
\end{array}
\right)  ^{2} & \cdots & \sum\limits_{m=0}^{L+M}c_{m}\left(
\begin{array}
[c]{ccc}%
M & m & L+M\\
0 & 0 & 0
\end{array}
\right)  ^{2}%
\end{array}
\right)
\end{equation}
and the matrix $B_{k}$\ is given by replacing the $k$-th column of the matrix
$A$ with
\end{paracol}

\[
\left[  -b_{0}\sum\limits_{m=0}^{L+M}c_{m}\left(
\begin{array}
[c]{ccc}%
0 & m & L+1\\
0 & 0 & 0
\end{array}
\right)  ^{2},-b_{0}\sum\limits_{m=0}^{L+M}c_{m}\left(
\begin{array}
[c]{ccc}%
0 & m & L+2\\
0 & 0 & 0
\end{array}
\right)  ^{2},\cdots,-b_{0}\sum\limits_{m=0}^{L+M}c_{m}\left(
\begin{array}
[c]{ccc}%
0 & m & L+M\\
0 & 0 & 0
\end{array}
\right)  ^{2}\right]  ^{T}:
\]%
\begin{equation}
B_{k}=\left(
\begin{array}
[c]{ccccc}%
\sum\limits_{m=0}^{L+M}c_{m}\left(
\begin{array}
[c]{ccc}%
1 & m & L+1\\
0 & 0 & 0
\end{array}
\right)  ^{2} & \cdots & -b_{0}\sum\limits_{m=0}^{L+M}c_{m}\left(
\begin{array}
[c]{ccc}%
0 & m & L+1\\
0 & 0 & 0
\end{array}
\right)  ^{2} & \cdots & \sum\limits_{m=0}^{L+M}c_{m}\left(
\begin{array}
[c]{ccc}%
M & m & L+1\\
0 & 0 & 0
\end{array}
\right)  ^{2}\\
\sum\limits_{m=0}^{L+M}c_{m}\left(
\begin{array}
[c]{ccc}%
0 & m & L+2\\
0 & 0 & 0
\end{array}
\right)  ^{2} & \cdots & -b_{0}\sum\limits_{m=0}^{L+M}c_{m}\left(
\begin{array}
[c]{ccc}%
0 & m & L+2\\
0 & 0 & 0
\end{array}
\right)  ^{2} & \cdots & \sum\limits_{m=0}^{L+M}c_{m}\left(
\begin{array}
[c]{ccc}%
M & m & L+2\\
0 & 0 & 0
\end{array}
\right)  ^{2}\\
\vdots & \vdots & \vdots & \vdots & \vdots\\
\sum\limits_{m=0}^{L+M}c_{m}\left(
\begin{array}
[c]{ccc}%
0 & m & L+M\\
0 & 0 & 0
\end{array}
\right)  ^{2} & \cdots & -b_{0}\sum\limits_{m=0}^{L+M}c_{m}\left(
\begin{array}
[c]{ccc}%
0 & m & L+M\\
0 & 0 & 0
\end{array}
\right)  ^{2} & \cdots & \sum\limits_{m=0}^{L+M}c_{m}\left(
\begin{array}
[c]{ccc}%
M & m & L+M\\
0 & 0 & 0
\end{array}
\right)  ^{2}%
\end{array}
\right)  .
\end{equation}

\begin{paracol}{2}
\switchcolumn
After obtaining the coefficient $b_{m}$, the coefficient $a_{n}$ is given by
Eq. (\ref{aneq}):%
\begin{equation}
a_{n}=2\sum_{l=0}^{M}\sum_{m=0}^{L+M}b_{l}c_{m}\left(
\begin{array}
[c]{ccc}%
l & m & n\\
0 & 0 & 0
\end{array}
\right)  ^{2},\text{ }0\leq n\leq L. \label{coefa}%
\end{equation}

\noindent\textit{Example:} Now we turn to the oscillation problem in the partial sum
(\ref{fNsin}) in section \ref{ex}.

\vspace*{-1cm}

\begin{wrapfigure}[20]{l}[0.95cm]{0.8\textwidth}
\includegraphics[width=0.8\textwidth]{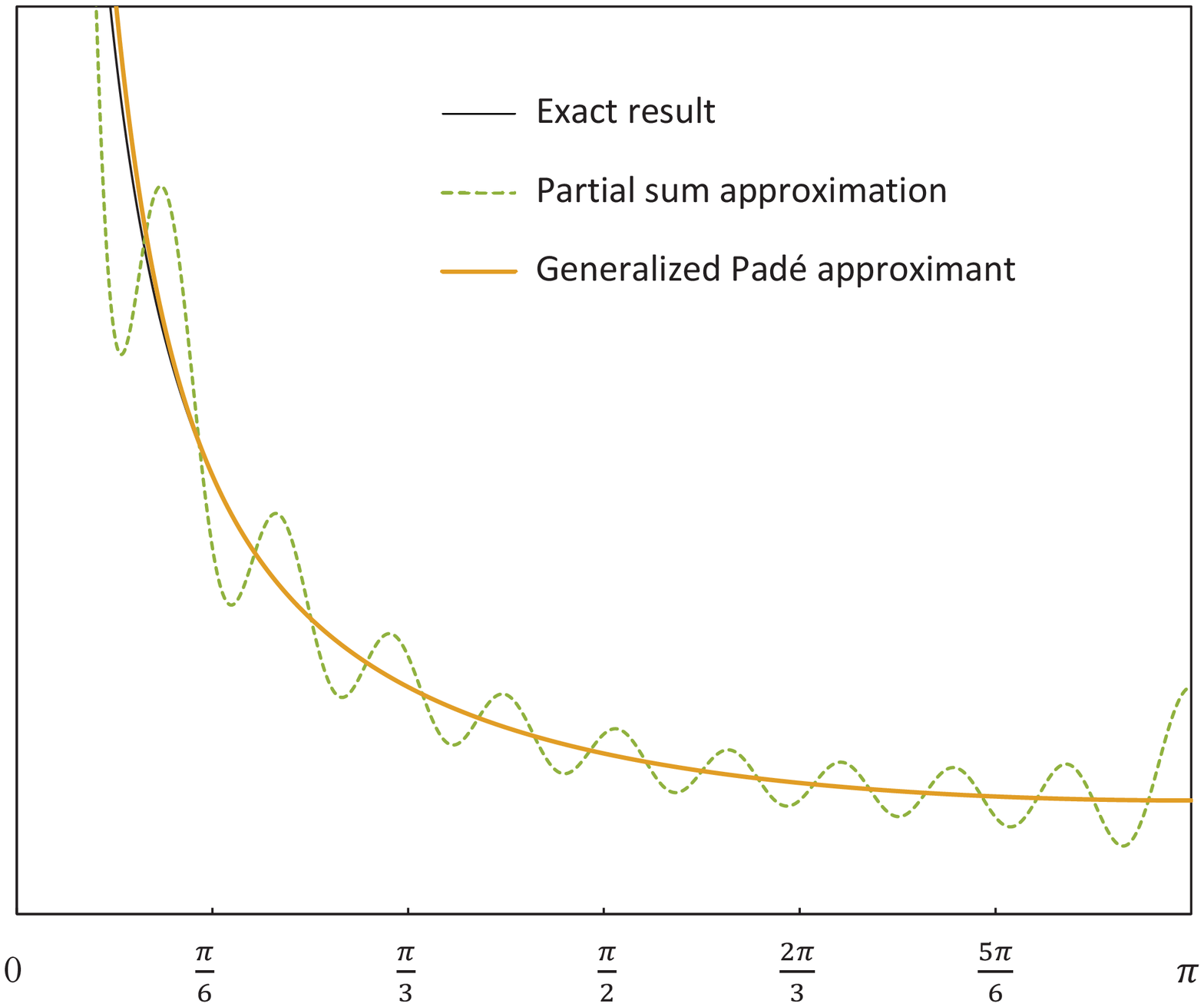}
\captionsetup{labelformat=empty} 
\caption{} 
\label{sumpade.eps} 
\end{wrapfigure}

\qquad\qquad\\
\qquad\qquad\\
\qquad\qquad\\
\qquad\qquad\\
\qquad\qquad\\
\qquad\qquad\\
\qquad\qquad\\
\qquad\qquad\\
\qquad\qquad\\
\qquad\qquad\\
\qquad\qquad\\
\qquad\qquad\\
\vspace*{4cm}
\qquad\qquad\\
\qquad\qquad\\
\qquad\qquad\\
\qquad\qquad\\
\qquad\qquad\\
\qquad\qquad\\
By replacing the polynomial in Eq. (\ref{fNsin}) with the rational expression
$f_{\left[  L/M\right]  }\left(  \theta\right)  $ defined by Eq.
(\ref{LPade}), we can greatly improve the oscillation problem.

Comparing with the exact solution, we can see that the rational expression
$f_{\left[  3/3\right]  }\left(  \theta\right)  $ constructed from the first
$6$ terms of the power series, i.e., taking $N=6$ in Eq. (\ref{fNsin}), gives
a very accurate result (see Fig. \ref{sumpade.eps}).
\switchcolumn
\vspace*{4.5cm}
{\color{Gray} 
\noindent\textbf{Figure \ref{sumpade.eps}.} Eliminating the incorrect oscillation in partial sum approximation of the
periodic function $f\left(  \theta\right)  =\frac{1}{2\sin\frac{\theta}{2}}$
by the generalized Pad{\'{e} }approximant. 
}
\switchcolumn

\section{Example: Cross section in quantum mechanics \label{QM}}

We take the scattering cross section in quantum mechanics as examples to
illustrate the method and its effectiveness.

\subsection{Coulomb potential: $V\left(  r\right)  =\alpha/r$}

Scattering on the coulomb potential has exact solutions.

The exact scattering amplitude of the Coulomb potential is
\cite{landau2013quantum}%
\begin{equation}
f\left(  \theta\right)  =-\frac{1}{2k^{2}\sin^{2}\frac{\theta}{2}}\frac
{\Gamma\left(  1+\frac{i}{k}\right)  }{\Gamma\left(  1-\frac{i}{k}\right)
}\exp\left(  -i\frac{2}{k}\ln\sin\frac{\theta}{2}\right)  . \label{fCouP}%
\end{equation}
Alternatively, the scattering amplitude of the Coulomb potential can also be
written as a sum as that in Eq. (\ref{fP}) \cite{landau2013quantum}:%

\begin{equation}
f\left(  \theta\right)  =\frac{1}{2ik}\sum_{l=0}^{\infty}\left(  2l+1\right)
\frac{\Gamma\left(  l+1+\frac{i}{k}\right)  }{\Gamma\left(  l+1-\frac{i}%
{k}\right)  }P_{l}\left(  \cos\theta\right)  .
\end{equation}
Though there is no oscillation in the exact result (\ref{fCouP}), an incorrect
oscillation appears in the partial sum.

As an example, corresponding to the first-$6$-term partial sum, we construct
the generalized Pad{\'{e}} approximant $f_{\left[  3/3\right]  }\left(
\theta\right)  $:
\begin{equation}
f_{\left[  3/3\right]  }\left(  \theta\right)  =\frac{P}{Q},
\end{equation}
where the numerator%
\begin{align}
P  &  =\sum_{n=0}^{3}a_{n}P_{n}\left(  \cos\theta\right) \nonumber\\
&  =\frac{\left(  \frac{4290272012250}{264699104689}+\frac{4320822501450}%
{264699104689}i\right)  \pi\operatorname{csch}\pi}{\Gamma(3-i)\Gamma
(4-i)}\nonumber\\
&  -\frac{\left(  \frac{86298291777600}{264699104689}-\frac{22045543159500}%
{264699104689}i\right)  \pi\operatorname{csch}\pi}{\Gamma(4-i)\Gamma
(5-i)}P_{1}\left(  \cos\theta\right) \nonumber\\
&  +\frac{2\Gamma\left(  1+i\right)  }{\Gamma\left(  2-i\right)  }\left(
\frac{4967588289069}{58498502136269}+\frac{7767606080115}{58498502136269}%
i\right)  P_{2}\left(  \cos\theta\right) \nonumber\\
&  -\frac{2\Gamma\left(  1+i\right)  }{\Gamma\left(  2-i\right)  }\left(
\frac{1846662362937}{166495736849381}-\frac{4075527063195}{166495736849381}%
i\right)  P_{3}\left(  \cos\theta\right)  ,
\end{align}
the denominator%
\begin{align}
Q  &  =\sum_{m=0}^{3}b_{n}P_{n}\left(  \cos\theta\right) \nonumber\\
&  =1-\left(  \frac{1570098416997}{1058796418756}+\frac{31193942505}%
{264699104689}i\right)  P_{1}\left(  \cos\theta\right) \nonumber\\
&  +2\left(  \frac{507748194515}{2117592837512}+\frac{66753972375}%
{1058796418756}i\right)  P_{2}\left(  \cos\theta\right) \nonumber\\
&  +2\left(  \frac{1746613473}{1058796418756}+\frac{384403371}{264699104689}%
i\right)  P_{3}\left(  \cos\theta\right)  .
\end{align}
and we take $\alpha=1$.

From Fig. (\ref{coulomb.eps})\ we can see that the incorrect oscillation has
been eliminated.

\subsection{Potential $V\left(  r\right)  =\alpha/r^{2}$}
\textit{Analytic result: Born approximation}. The scattering amplitude in
quantum mechanical scattering is \cite{landau2013quantum}
\begin{equation}
f\left(  \theta\right)  =\frac{1}{2ik}\sum_{l=0}^{\infty}\left(  2l+1\right)
\left(  e^{2i\delta_{l}}-1\right)  P_{l}\left(  \cos\theta\right)  .
\end{equation}
Under the small phase shift approximation, the scattering amplitude can be
written as%
\begin{equation}
f\left(  \theta\right)  =\frac{1}{k}\sum_{l=0}^{\infty}\left(  2l+1\right)
\delta_{l}P_{l}\left(  \cos\theta\right)  \label{fthsmall}%
\end{equation}

By the Born approximation, the partial wave scattering phase shift is
\begin{equation}
\delta_{l}^{\text{Born}}=-k\int_{0}^{\infty}j_{l}^{2}\left(  kr\right)
V\left(  r\right)  r^{2}dr, \label{deltaborn}%
\end{equation}
where $V\left(  r\right)  $ is the potential and $j_{l}\left(  z\right)  $ is
the spherical Bessel functions.

\vspace*{-1cm}
\begin{wrapfigure}[8]{l}[0.95cm]{0.8\textwidth}
\includegraphics[width=0.8\textwidth]{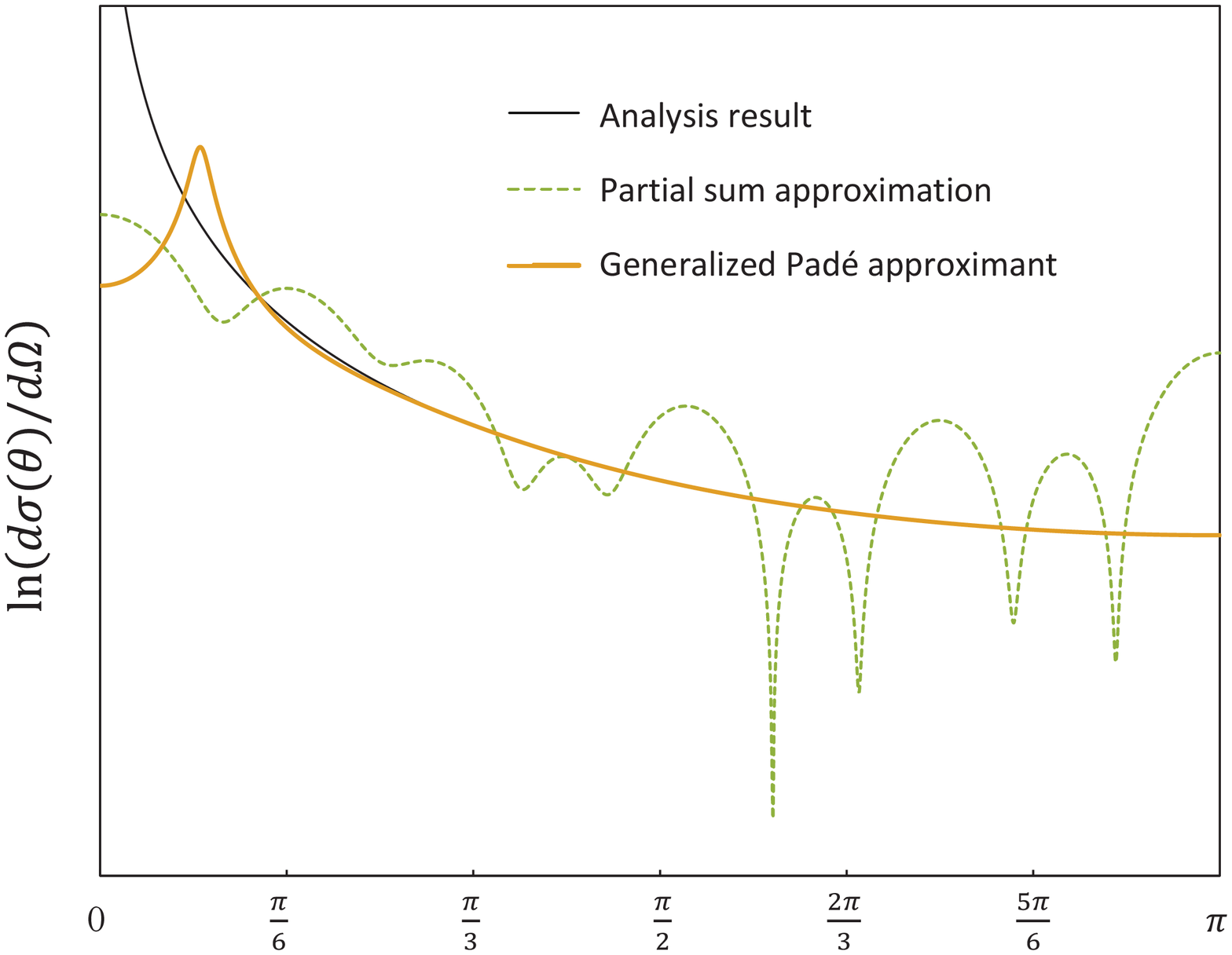}
\captionsetup{labelformat=empty} 
\caption{} 
\label{coulomb.eps} 
\end{wrapfigure}
\qquad\qquad
\qquad\qquad
\qquad\qquad
\qquad\qquad
\\ 
\qquad\qquad
\\
\qquad\qquad
\\
\vspace*{8cm}
\qquad\qquad
\\
\qquad\qquad
\\
\qquad\qquad
\\\qquad\qquad
\\

\switchcolumn
\newpage
\qquad
\newpage
\vspace*{8.5cm}
{\color{Gray} 
\noindent\textbf{Figure \ref{coulomb.eps}.} Differential scattering cross sections of the Coulomb potential. 
}
\switchcolumn

\vspace*{1cm}
The Born approximation scattering amplitude, by substituting Eq.
(\ref{deltaborn}) into Eq. (\ref{fthsmall}), is
\begin{equation}
f^{\text{Born}}\left(  \theta\right)  =-\int_{0}^{\infty}drr^{2}V\left(
r\right)  \sum_{l=0}^{\infty}\left(  2l+1\right)  j_{l}^{2}\left(  kr\right)
P_{l}\left(  \cos\theta\right)  . \label{fBornsum}%
\end{equation}
In this first order approximation, we encounter the sum of the form (\ref{fP}).

The sum here, however, can be performed exactly,%
\begin{equation}
\sum_{l=0}^{\infty}\left(  2l+1\right)  j_{l}^{2}\left(  kr\right)
P_{l}\left(  \cos\theta\right)  =\frac{\sin\left(  qr\right)  }{qr}%
\end{equation}
with $q=2k\sin\frac{\theta}{2}$, so the amplitude (\ref{fBornsum}) becomes%
\begin{equation}
f^{\text{Born}}\left(  \theta\right)  =-\int_{0}^{\infty}r^{2}drV\left(
r\right)  \frac{\sin\left(  qr\right)  }{qr}. \label{fBorn}%
\end{equation}
This enables us to check the validity of the method through directly comparing
the approximate result given by the partial sum with the analytic result.

For the potential
\begin{equation}
V\left(  r\right)  =\frac{\alpha}{r^{2}}, \label{potentialr2}%
\end{equation}
by Eq. (\ref{fBorn}) we have
\begin{equation}
f^{\text{Born}}\left(  \theta\right)  =-\frac{\pi\alpha}{4k\sin\frac{\theta
}{2}}. \label{fBornAn}%
\end{equation}

\textit{Partial sum approximation:} The partial sum approximation is%
\begin{equation}
f_{N}\left(  \theta\right)  \simeq\frac{1}{k}\sum_{l=0}^{N}\left(
2l+1\right)  \delta_{l}^{\text{Born}}P_{l}\left(  \cos\theta\right)  .
\label{fN}%
\end{equation}

This partial sum is of the form of Eq. (\ref{Nterms}), which leads to an
incorrect oscillation as in the example given in section \ref{ex}. This
oscillation does not appear in the result given by (\ref{fBornAn}). Such an
oscillation cannot be eliminated by keeping more terms.

\textit{Eliminating oscillation.} We now use the generalized Pad{\'{e}}
approximant constructed in section \ref{Scheme} to eliminate the incorrect oscillation.

\vspace*{-1cm}

\begin{wrapfigure}[20]{l}[0.95cm]{0.8\textwidth}
\includegraphics[width=0.8\textwidth]{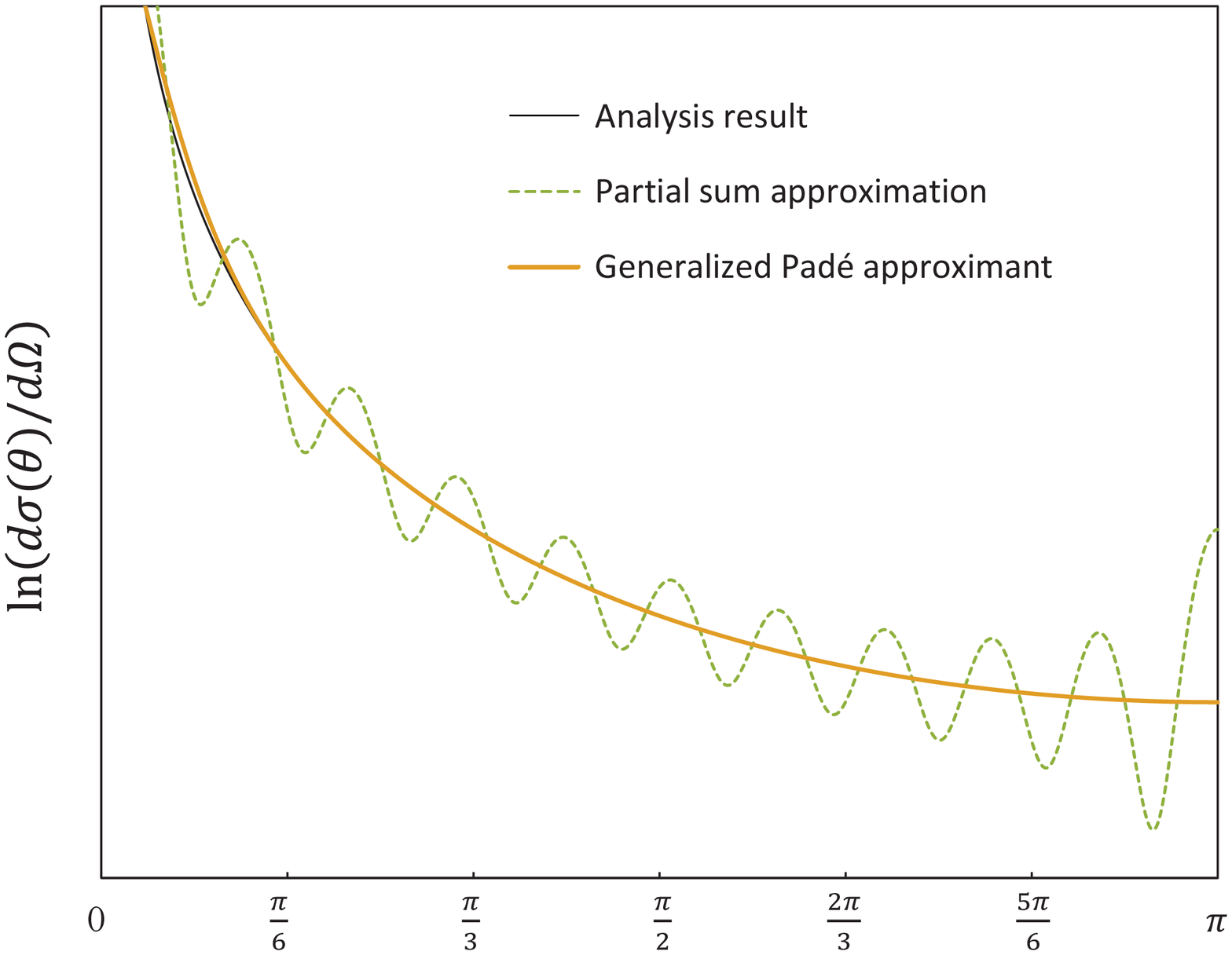}
\captionsetup{labelformat=empty} 
\caption{} 
\label{n2Born.eps} 
\end{wrapfigure}
\qquad\qquad
\qquad\qquad
\qquad\qquad
\qquad\qquad
\\ 
\qquad\qquad
\\
\qquad\qquad
\\
\qquad\qquad
\\
\qquad\qquad
\\
\qquad\qquad
\\\qquad\qquad
\\
\\\qquad\qquad
\\
\\\qquad\qquad
\\
\\\qquad\qquad
\\
\\\qquad\qquad
\\
\\\qquad\qquad
\\\qquad\qquad
\\

\newpage
Taking the first $6$ terms, i.e., $N=6$, in Eq. (\ref{fN}) as an example,
\begin{equation}
f_{N=6}\left(  \theta\right)  =\frac{1}{k}\sum_{l=0}^{6}\left(  2l+1\right)
\delta_{l}^{\text{Born}}P_{l}\left(  \cos\theta\right)
\end{equation}
\switchcolumn
\newpage
\vspace*{13.5cm}
{\color{Gray} 
\noindent\textbf{Figure \ref{n2Born.eps}.} Differential scattering cross sections of the potential $V\left(  r\right)
=\frac{\alpha}{r^{2}}$. 
}
\switchcolumn

The generalized Pad{\'{e}} approximant (\ref{PadeL}) is%
\begin{equation}
f_{\left[  3/3\right]  }\left(  \theta\right)  =\frac{-\frac{184224\pi
}{5948545}+\frac{988256\pi}{29742725}P_{1}\left(  \cos\theta\right)
-\frac{257088\pi}{65433995}P_{2}\left(  \cos\theta\right)  -\frac{1036608\pi
}{4253209675}P_{3}\left(  \cos\theta\right)  }{1-\frac{16158513}%
{11897090}P_{1}\left(  \cos\theta\right)  +\frac{854777}{2379418}P_{2}\left(
\cos\theta\right)  +\frac{11424}{5948545}P_{3}\left(  \cos\theta\right)  },
\end{equation}
where we take $\alpha=1$.

It can be seen from Fig. (\ref{n2Born.eps}) that the generalized Pad{\'{e}}
approximant works well.

\section{Example: Scattering in Reissner-Nordstr\"{o}m spacetime \label{RN}}

In the calculation of the scattering cross section in the
Reissner-Nordstr\"{o}m spacetime, we also encounter the incorrect oscillation
of the partial sum. Such an oscillation can also be eliminated by the
generalized Pad{\'{e} }approximant.
\end{paracol}
\begin{figure}[h]
\centering\includegraphics[width=1.05\textwidth]{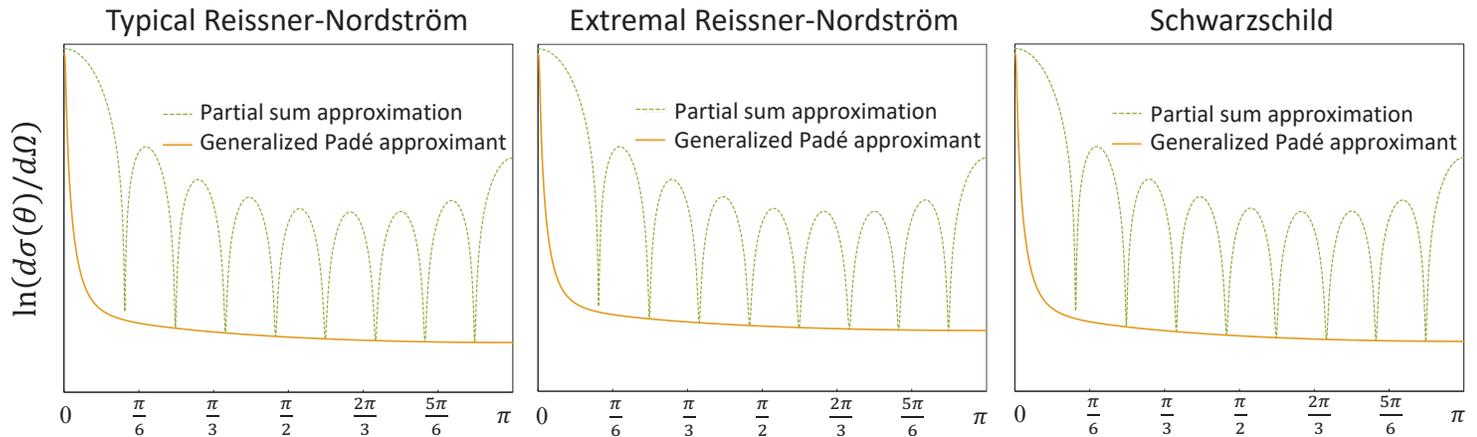}\caption{{}The oscillation
in the partial sum approximation and the elimination of the oscillation: the
typical Reissner-Nordstr\"{o}m case $Q/M=1/2$, the extremal
Reissner-Nordstr\"{o}m case $Q/M\simeq1$, and the Schwarzschild case
$Q/M\simeq0$.}%
\label{3RN.eps}%
\end{figure}
\begin{paracol}{2}
\switchcolumn
\bigskip
For scattering cross section in the Reissner-Nordstr\"{o}m spacetime, the
scattering amplitude is \cite{pike2008scatteringPage}%
\begin{equation}
f\left(  \theta\right)  =\frac{1}{2i\omega}\sum_{l=0}^{\infty}\left(
2l+1\right)  \left(  e^{2i\delta_{l}}-1\right)  P_{l}\left(  \cos
\theta\right)  ,\label{fTheta}%
\end{equation}
and the differential scattering cross section is
\begin{equation}
\sigma\left(  \theta\right)  =\left\vert f\left(  \theta\right)  \right\vert
^{2}.
\end{equation}
Here the zeroth-order phase shift is \cite{li2021scalar}
\begin{align}
\delta_{l}^{\left(  0\right)  } &  =\frac{l\pi}{2}+\frac{r_{+}+r_{-}}{2}%
\eta-\left(  r_{+}+r_{-}\right)  \eta\ln2\nonumber\\
&  =\frac{l\pi}{2}+M\eta-2M\eta\ln2+2M\eta\ln\frac{\sqrt{M^{2}-Q^{2}}}%
{M},\label{ps0}%
\end{align}
and the first-order phase shift is
\begin{align}
\delta_{l}^{\left(  1\right)  } &  =-\arctan\frac{\displaystyle\frac{1}{\eta
}\int_{r_{+}}^{\infty}\sin^{2}\left(  \eta r_{\ast}\right)  \frac{dr_{\ast}%
}{dr}V_{l}^{\text{eff}}dr}{1+\displaystyle\frac{1}{\eta}\int_{r_{+}}^{\infty
}\sin\left(  2\eta r_{\ast}\right)  \frac{dr_{\ast}}{dr}V_{l}^{\text{eff}}%
dr}+\left(  r_{+}+r_{-}\right)  \eta\ln\frac{r_{+}-r_{-}}{r_{+}+r_{-}%
}\nonumber\\
&  =-\arctan\frac{\displaystyle\frac{1}{\eta}\int_{r_{+}}^{\infty}\sin
^{2}\left(  \eta r_{\ast}\right)  \frac{dr_{\ast}}{dr}V_{l}^{\text{eff}}%
dr}{1+\displaystyle\frac{1}{\eta}\int_{r_{+}}^{\infty}\sin\left(  2\eta
r_{\ast}\right)  \frac{dr_{\ast}}{dr}V_{l}^{\text{eff}}dr}+\left(  r_{+}%
+r_{-}\right)  \eta\ln\frac{r_{+}-r_{-}}{r_{+}+r_{-}}.\label{ps1}%
\end{align}
Here $r_{\pm}=M\pm\sqrt{M^{2}-Q^{2}}$ are the horizons of the
Reissner-Nordstr\"{o}m spacetime, $r_{\ast}=r+\frac{r_{+}^{2}}{r_{+}-r_{-}}%
\ln\left(  \frac{r}{r_{+}}-1\right)  -\frac{r_{-}^{2}}{r_{+}-r_{-}}\ln\left(
\frac{r}{r_{-}}-1\right)  \ $is the tortoise coordinate of the
Reissner-Nordstr\"{o}m spacetime, $\eta^{2}=\omega^{2}-\mu^{2}$ with
$\omega^{2}$ the energy of the incident particle and $\mu$\ the mass of the
particle, and the effective potential $V_{l}^{\text{eff}}=(1-\frac{r_{+}}%
{r})(1-\frac{r_{-}}{r})[\frac{l\left(  l+1\right)  }{r^{2}}+(\frac{r_{+}%
+r_{-}}{r^{3}}-\frac{2r_{+}r_{-}}{r^{4}})]+\mu^{2}(\frac{r_{+}r_{-}}{r^{2}%
}-\frac{r_{+}+r_{-}}{r})$.

The scattering amplitude by the phase shift given by\ substituting Eqs.
(\ref{ps0}) and (\ref{ps1}) into Eq. (\ref{fTheta}), up to $l=6$, is%
\begin{equation}
f_{N=6}\left(  \theta\right)  =\frac{1}{2i\omega}\sum_{l=0}^{6}\left(
2l+1\right)  \exp2i\left(  \delta_{l}^{\left(  0\right)  }+\delta_{l}^{\left(
1\right)  }\right)  P_{l}\left(  \cos\theta\right)  .\label{f6RN}%
\end{equation}
In this partial sum, an incorrect oscillation appears, see Fig. (\ref{3RN.eps}).

In order to eliminate the oscillation, instead of the polynomial approximation
(\ref{f6RN}), we construct the generalized Pad{\'{e} }approximant as follows:%
\begin{equation}
f_{\left[  3/3\right]  }\left(  \theta\right)  =\frac{\sum_{n=0}^{3}a_{n}%
P_{n}\left(  \cos\theta\right)  }{\sum_{m=0}^{3}b_{m}P_{m}\left(  \cos
\theta\right)  }.
\end{equation}

Concretely, as examples, for parameters $\eta=10^{-4}$, $\mu=10^{-6}$, and
$M=10$, we have the following generalized Pad{\'{e} }approximants.

For the typical Reissner-Nordstr\"{o}m case $Q/M=1/2$,%
\end{paracol}
\begin{align}
f_{\left[  3/3\right]  }\left(  \theta\right)   &
=[(1156.89-156.71i)-(1444.37-185.02i)P_{1}\left(  \cos\theta\right)
+(288.3-28.64i)P_{2}\left(  \cos\theta\right)  \nonumber\\
&  +(5.70-1.40i)P_{3}\left(  \cos\theta\right)  ]/\left[  1-1.33P_{1}\left(
\cos\theta\right)  +(0.33-0.001i)P_{2}\left(  \cos\theta\right)  \right]  ,
\end{align}
\begin{paracol}{2}
\switchcolumn
for the extremal Reissner-Nordstr\"{o}m case $Q/M=0.99$,%
\end{paracol}
\begin{align}
f_{\left[  3/3\right]  }\left(  \theta\right)   &
=[(1422.85-266.46i)-(1796.91-323.3i)P_{1}\left(  \cos\theta\right)
+(374.62-57.17i)P_{2}\left(  \cos\theta\right)  \nonumber\\
&  +(5.77-1.9i)P_{3}\left(  \cos\theta\right)  ]/[1-(1.33-0.001i)P_{1}\left(
\cos\theta\right)  +(0.33-0.001i)P_{2}\left(  \cos\theta\right)  ],
\end{align}
\begin{paracol}{2}
\switchcolumn
and for the Schwarzschild case $Q/M=10^{-4}$,%
\end{paracol}
\begin{align}
f_{\left[  3/3\right]  }\left(  \theta\right)   &
=[(1165.19-157.75i)-(1453.4-185.79i)P_{1}\left(  \cos\theta\right)
+(289.06-28.4i)P_{2}\left(  \cos\theta\right)  \nonumber\\
&  +(5.83-1.43i)P_{3}\left(  \cos\theta\right)  ]/[1-(1.33-0.001i)P_{1}\left(
\cos\theta\right)  +(0.33-0.001i)P_{2}\left(  \cos\theta\right)  ].
\end{align}

\begin{paracol}{2}
\switchcolumn
\section{Conclusion \label{Conclusion}}

The approach suggested in this appendix can be used for eliminating the
oscillation in the truncated Fourier series, i.e., the partial sum
approximation of the Fourier series. Expanding a periodic function needs a
complete set consisting of the periodic function basis. The Fourier series
chooses the sine and cosine functions as the basis, and the spherically
symmetric scattering chooses the Legendre polynomial with the variable
$\cos\theta$ as the basis. The complete set of the Legendre polynomial with
the variable $\cos\theta$ is a rearrangement of the complete set of the sine
and cosine functions --- the Fourier case. Therefore, the Gibbs phenomenon of
the Fourier series will be transferred to the series with the Legendre
polynomial basis. The key in constructing the rational approximant for
eliminating the oscillation in the Gibbs phenomenon is that the rational
approximant should be constructed from the basis of the corresponding series.
For example, the rational approximant corresponding to the Fourier series
should be constructed by the sine and cosine functions, and the rational
approximant corresponding to the series with the Legendre polynomial should be
constructed by the Legendre polynomial, and so on. Only in this way can the
efficiency of the calculation be guaranteed.

\acknowledgments

We are very indebted to Dr G. Zeitrauman for his encouragement. This work is supported in part by Special Funds for theoretical physics Research Program of the NSFC under Grant No.
11947124, and NSFC under Grant Nos. 11575125 and 11675119.

\nolinenumbers

\reftitle{References}




\providecommand{\href}[2]{#2}\begingroup\raggedright\endgroup


\begin{thebibliography}{1}

\bibitem{tian2021pade}
Y.-H. Tian, W.-D. Li, Y.~Shen, and W.-S. Dai, {\it Pad{\'e} approximant
  approach to singular properties of quantum gases: The ideal cases},  {\em
  Communications in Theoretical Physics} (2021).

\bibitem{olver2010nist}
F.~W. Olver, D.~W. Lozier, R.~F. Boisvert, and C.~W. Clark, {\em NIST handbook
  of mathematical functions}.
\newblock Cambridge University Press, 2010.

\bibitem{landau2013quantum}
L.~Landau and E.~Lifshitz, {\em Quantum Mechanics: Non-Relativistic Theory}.
\newblock Teoreticheskai︠a︡ fizika.

\bibitem{pike2008scatteringPage}
E.~Pike and P.~Sabatier, {\em Scattering: Scattering and Inverse Scattering in
  Pure and Applied Science}.
\newblock Academic, "2008 \,"" and Page 1613".

\bibitem{li2021scalar}
S.-L. Li, Y.-Y. Liu, W.-D. Li, and W.-S. Dai, {\it Scalar field in
  Reissner-Nordstr{\"o}m spacetime: Bound state and scattering state (with
  appendix on eliminating oscillation in partial sum approximation of periodic
  function)},  {\em Annals of Physics} (2021) 168578.

\end{thebibliography}
\end{paracol}
%


\end{document}